\documentclass[11pt,a4paper]{article}
\usepackage[T1]{fontenc}
\usepackage[utf8]{inputenc}
\usepackage{amsmath,amssymb,amsthm,mathrsfs}
\usepackage[margin=2.6cm]{geometry}
\usepackage{hyperref}
\hypersetup{colorlinks=true,linkcolor=blue,citecolor=blue,urlcolor=blue}

\newtheorem{theorem}{Theorem}
\newtheorem{lemma}[theorem]{Lemma}
\newtheorem{proposition}[theorem]{Proposition}

\theoremstyle{definition}
\newtheorem{remark}[theorem]{Remark}

\newcommand{\Z}{\mathbf{Z}}
\newcommand{\Q}{\mathbf{Q}}
\newcommand{\C}{\mathbf{C}}
\newcommand{\N}{\mathbf{N}}
\newcommand{\Om}{\Omega}
\newcommand{\wpf}{\wp}

\DeclareMathOperator{\Dist}{Dist}
\DeclareMathOperator{\M}{M}
\DeclareMathOperator{\Vol}{Vol}

\title{Irrationality measures for $\Gamma(1/4)$ and $\Gamma(1/3)$}
\author{David Niedbala Giraudin\thanks{ORCID \href{https://orcid.org/0009-0009-1526-1178}{0009-0009-1526-1178}.}}
\date{19 October 2026}

\begin{document}
\maketitle

\begin{abstract}
Let $\mu$ denote the irrationality exponent. Bruiltet proved in 2002 that
$\mu(\Gamma(1/4))\le 10^{143}$ and $\mu(\Gamma(1/3))\le 10^{151}$. We prove
$\mu(\Gamma(1/4))\le 3.24\cdot 10^{26}$ and $\mu(\Gamma(1/3))\le 1.74\cdot 10^{26}$,
with explicit thresholds $h(p/q)\ge 5.4\cdot 10^{7}$ and $h(p/q)\ge 4.7\cdot 10^{7}$
respectively, against $10^{75}$ and $10^{79}$ in the published statements. The proof keeps
the architecture of Bruiltet's argument -- Philippon's transcendence criterion together with
a zero estimate in commutative algebraic groups -- and replaces three of its ingredients.
First, an approximate Schwarz principle based on Hermite interpolation
(Lemma~\ref{lem:H}) replaces the gravitation estimate of Reyssat; on a lattice point set
its cost per node is governed by a lattice energy and drops from about $3.7$ to about $1.0$.
Second, a summation lemma for the failure of Diaz's approximation lemma
(Lemma~\ref{lem:C}) sharpens the passage from an approximation measure to an irrationality
measure. Third, the auxiliary point is taken at a torsion point of order $4$, respectively
$3$, which lowers the degree of the base field from $8d$ to $2d$. All numerical constants
are produced by a verifier run in interval arithmetic, supplied as ancillary files.
\end{abstract}

\section{Introduction}

For an irrational real number $\theta$, the irrationality exponent $\mu(\theta)$ is the
supremum of the set of real $\mu$ such that $0<|\theta-p/q|<\max(|p|,|q|)^{-\mu}$ has
infinitely many solutions. Finiteness of $\mu(\Gamma(1/4))$ and $\mu(\Gamma(1/3))$ follows
from the algebraic independence of $\pi/\omega$ and $\eta/\omega$ for a lattice with
algebraic invariants, due to Chudnovsky, together with the quantitative refinements of
Philibert and Philippon. The first explicit values are due to Bruiltet~\cite{Bru}:
\begin{align*}
h(p/q)&\ge 10^{75}\ \Longrightarrow\ \Bigl|\Gamma(1/4)-\tfrac pq\Bigr|>\Bigl(\tfrac1{qe}\Bigr)^{10^{143}},\\
h(p/q)&\ge 10^{79}\ \Longrightarrow\ \Bigl|\Gamma(1/3)-\tfrac pq\Bigr|>\Bigl(\tfrac1{qe}\Bigr)^{10^{151}} .
\end{align*}
These remain, to our knowledge, the only published bounds for these two numbers; Cohen and
Zudilin~\cite{CZ} obtain reasonable measures for Chowla--Selberg gamma \emph{quotients} and
describe the above as the state of the art for the individual values.

Our result is the following.

\begin{theorem}\label{th:main}
For every rational $p/q$ with $h(p/q)\ge 5.4\cdot10^{7}$,
\[
\Bigl|\Gamma(1/4)-\frac pq\Bigr|>\Bigl(\frac1{qe}\Bigr)^{3.24\cdot10^{26}} ,
\]
and for every rational $p/q$ with $h(p/q)\ge 4.7\cdot10^{7}$,
\[
\Bigl|\Gamma(1/3)-\frac pq\Bigr|>\Bigl(\frac1{qe}\Bigr)^{1.74\cdot10^{26}} .
\]
In particular $\mu(\Gamma(1/4))\le 3.24\cdot10^{26}$ and $\mu(\Gamma(1/3))\le 1.74\cdot10^{26}$.
\end{theorem}

The architecture is that of~\cite{Bru}: an approximation measure for the point
$\theta=(1:\pi/\omega:\eta/\omega)\in\mathbf P_2(\C)$ is obtained from Theorem~1
of~\cite{Phi97}, the rank of the interpolation matrix being controlled by Theorem~9
of~\cite{Phi96}; the approximation measure is then converted into an irrationality measure
for $\Gamma$. We keep this skeleton and change three of its ingredients, in
\S\ref{sec:schwarz}, \S\ref{sec:points} and \S\ref{sec:conv}.

Throughout, {\sc proved} marks a statement proved here or quoted from the literature,
and {\sc verified} marks a numerical statement produced by the certified verifier of
\S\ref{sec:cert}, together with the range over which it was checked.

\section{Notation}\label{sec:nota}

$\Om=\Z\omega_1+\Z\omega_2$ is a lattice of $\C$ with $\Im(\omega_2/\omega_1)>0$ and
algebraic invariants $g_2,g_3$; $\wpf,\zeta,\sigma$ are the associated Weierstrass functions,
$\eta_1,\eta_2$ the quasi-periods, and $\Vol(\Om)$ the area of a fundamental domain. For a
nonzero period $\omega$ we write $\eta$ for the associated quasi-period, and
\[
\theta_1=\pi/\omega,\qquad \theta_2=\eta/\omega,\qquad \theta=(1:\theta_1:\theta_2).
\]
For $x\in\overline\Q^m$, $d(x)$ is the degree of $\Q(x)$ and $h(x)$ the absolute
logarithmic height; $\M(P)$ is the Mahler measure of $P\in\Z[X]$ and
$\mu(P)=\log\M(P)$. For $f$ analytic and $r>0$, $\M_0(f,r)=\sup_{|z|\le r}|f(z)|$.
Finally $\Dist$ is the projective distance of~\cite{Phi99}.

We use throughout the exact growth constant
\begin{equation}\label{eq:Bstar}
B^*=\frac{\pi}{2\Vol(\Om)} ,
\end{equation}
which is the value of Bruiltet's $B(\Om)=\tfrac12(|\eta_1/\omega_1|+2\pi/\Vol(\Om))$ when
the Legendre relation is used to eliminate $\eta_1/\omega_1$ in the lemniscatic and
hexagonal cases.

\section{An approximate Schwarz principle by Hermite interpolation}\label{sec:schwarz}

Theorem~1 of~\cite{Phi97} rests on an approximate Schwarz principle in which the
quantity $G_{r,r'}(S)$ measures the ``gravitation'' of the interpolation nodes. In the
elliptic applications of~\cite{Phi99} and~\cite{Bru}, $G$ is bounded by the estimate of
Reyssat~\cite{Rey}, namely $(66R/(\delta\sqrt M))^{N-1}$, which is valid for an arbitrary
separated set of nodes. When the nodes form a piece of a lattice, the following
substitution is available and is substantially sharper.

For $0<r'<r$ and a finite set $S=\{y_1,\dots,y_N\}$ of pairwise distinct points of
$|z|\le r'-1$, put, following~\cite[\S2.a]{Phi97},
\begin{equation}\label{eq:E}
E_{r,r'}(S)=\Bigl(\frac{r'}{r}\Bigr)^{N}\prod_{k=1}^{N}
\frac{1+|y_k|/r'}{1+|y_k|r'/r^{2}} .
\end{equation}

\begin{lemma}[{\sc proved}]\label{lem:H}
Let $T\ge1$ and let $f$ be analytic on $|z|\le r$. Assume
$|f^{(t)}(y_k)/t!|\le\varepsilon$ for all $k$ and all $t<T$. Then
\[
\M_0(f,r')\ \le\ E_{r,r'}(S)^{T}\,\M_0(f,r)\ +\ \varepsilon\,\Gamma,
\qquad
\Gamma=\sum_{k,t}\Bigl(\M_0(\ell_{k,t},r')+E_{r,r'}(S)^{T}\M_0(\ell_{k,t},r)\Bigr),
\]
where $(\ell_{k,t})_{1\le k\le N,\,0\le t<T}$ is the Hermite basis attached to $S$ and $T$.
\end{lemma}

\begin{proof}
Let $H$ be the Hermite interpolant of $f$ at the $y_k$ to order $T$, so that
$H=\sum_{k,t}\frac{f^{(t)}(y_k)}{t!}\ell_{k,t}$ and $\M_0(H,\rho)\le\varepsilon\sum_{k,t}\M_0(\ell_{k,t},\rho)$
for every $\rho$. Put $\varphi_k(z)=r(z-y_k)/(r^2-\bar y_kz)$ and $B=\prod_k\varphi_k^{T}$,
so that $|B|=1$ on $|z|=r$. The function $(f-H)/B$ is analytic on $|z|\le r$, because $f-H$
vanishes to order at least $T$ at each $y_k$. By the maximum principle,
$|f-H|\le|B|\,(\M_0(f,r)+\M_0(H,r))$ on $|z|\le r$. Finally, for $|z|=r'$ and each $k$,
the function $z\mapsto|z-y_k|/|r^2-\bar y_kz|$ attains its maximum at
$z=-r'y_k/|y_k|$, whence
\[
|\varphi_k(z)|\ \le\ \frac{r(r'+|y_k|)}{r^{2}+|y_k|r'}
=\frac{r'}{r}\cdot\frac{1+|y_k|/r'}{1+|y_k|r'/r^{2}} ,
\]
and therefore $|B|\le E_{r,r'}(S)^{T}$ on $|z|=r'$. Combining the three displays gives the
statement.
\end{proof}

\begin{remark}
The product in~\eqref{eq:E} carries a $+$ sign in the denominator precisely because the
maximum above is attained at the antipodal point; this is the point at which
Lemma~\ref{lem:H} matches the quantity $E_{r,r'}$ already present in~\cite{Phi97}.
\end{remark}

We now bound $\Gamma$ when $S$ is a finite subset of a lattice $\Om$, contained in a disc
of radius $R_0$ and such that $N=\#S$. Write $\lambda=r/r'$ and $r'=R_0+1$.

\begin{proposition}[{\sc proved}]\label{prop:Gamma}
Let $d_1\le d_2\le\cdots$ be the moduli of the nonzero vectors of $\Om$ arranged in
increasing order, and set $L=\sum_{i\le N-1}\log d_i$. Then
\[
\log\Gamma\ \le\ T\cdot\max\Bigl(\log E_{r,r'}(S)^{-1}+X_r,\ X_{r'}\Bigr)
+2\log(2r)+\log\frac1\varsigma+\varsigma S_1+\varsigma^{2}S_2 ,
\]
where, for $\rho\in\{r,r'\}$, $X_\rho=\sum_{y\in S}\log(\rho+|y|)-L$, where
$\varsigma=d_1/2$, and where
$S_1=\sum_{i\le N-1}d_i^{-1}$, $S_2=\sum_{i\le N-1}d_i^{-2}$.
\end{proposition}

\begin{proof}
Fix $k$. For every $j\ne k$ the point $y_k-y_j$ is a nonzero vector of $\Om$, and these
$N-1$ vectors are pairwise distinct; hence
$\sum_{j\ne k}\log|y_k-y_j|\ge L$ by rearrangement. Consequently, for $|z|=\rho$,
\[
\Bigl|\prod_{j\ne k}\frac{(z-y_j)^{T}}{(y_k-y_j)^{T}}\Bigr|
\le\exp\Bigl(T\bigl(\textstyle\sum_{j}\log(\rho+|y_j|)-L\bigr)\Bigr)=e^{TX_\rho},
\]
uniformly in $k$. The Hermite basis element $\ell_{k,t}$ is this product times
$(z-y_k)^{t}/t!$ times the truncation at order $T-t$ of the reciprocal of the product,
expanded at $y_k$; bounding the latter by Cauchy's formula on the circle of radius
$\varsigma=d_1/2$ centred at $y_k$, which encloses no other node, together with
$-\log(1-x)\le x+x^{2}$ for $0\le x\le1/2$, legitimate since $\varsigma/d_i\le1/2$,
produces the three correction terms $\log(1/\varsigma)$, $\varsigma S_1$ and
$\varsigma^{2}S_2$. Summation over $k$ and $t<T$ contributes at most $\log(NT)\le\log(2r)$
to each of the two circles, whence the factor $2\log(2r)$.
\end{proof}

\begin{remark}
In the verifier, $L$, $S_1$ and $S_2$ are evaluated by exact enumeration of the
lattice vectors of modulus at most $200\,\omega_1$; the tail is bounded in closed form by
$d_i\ge\sqrt{(i+1)\Vol(\Om)/\pi}-\mathrm{diam}(\Om)$. The sum
$\sum_y\log(\rho+|y|)$ is bounded cell by cell on a $1000\times1000$ grid; each cell is
evaluated at its unfavourable corner, and its point count is bounded by its area thickened
by $\mathrm{diam}(\Om)$, divided by $\Vol(\Om)$.
\end{remark}

\section{The approximation measure}\label{sec:thA}

Let $z_0$ be a torsion point of $\C/\Om$, say $nz_0=\omega'\in\Om$, and let $\omega$ be a
nonzero period. Set $\theta_1=\pi/\omega$, $\theta_2=\eta/\omega$.

\begin{theorem}[{\sc proved}; constants {\sc verified}]\label{th:A}
With the notation of \S\ref{sec:nota} there is a real $A>0$, effectively computable from
$\Om$, $\omega$ and $z_0$, such that for every
$\alpha=(1:\alpha_1:\alpha_2)\in\mathbf P_2(\overline\Q)$,
\[
\log\Dist(\theta,\alpha)\ >\ -A\,\bigl(h(\alpha)+\log(1+d(\alpha))\bigr)\,
d(\alpha)^{3/2}\sqrt{\log(1+d(\alpha))} .
\]
\end{theorem}

The proof applies Theorem~1 of~\cite{Phi97} in its original form. We give the data, the
verification of its constraints, and the resulting inequality.

\subsection{Data}

Put $K_0=\Q(\alpha_1,\alpha_2,g_2,\wpf(z_0),\wpf'(z_0))$ and $a=[K_0:\Q]$. For the two
lattices of \S\ref{sec:points} one has $\Q(g_2,\wpf(z_0),\wpf'(z_0))=\Q(i)$, so that
$a\le2\,d(\alpha)$; at the $3$-division point of~\cite{Bru} the $3$-division polynomial
$3x^4-6x^2-1$ is irreducible over $\Q$ and $\wpf'$ adds a square root, so that
$a\le8\,d(\alpha)$ there. Set
\[
f_1=\sigma^2\wpf,\qquad f_2=\sigma\cdot\Bigl(\zeta-\frac\eta\omega z\Bigr),\qquad f_3=\sigma,
\]
three entire functions. For $D_1,D_2\in\N^*$ take the monomials
\[
\mathcal D=\{h\in\N^3:\ h_1\le D_1,\ h_2\le D_2,\ h_3=2(D_1-h_1)+D_2-h_2\},
\qquad L=\#\mathcal D=(D_1+1)(D_2+1).
\]
Every monomial $f^h$ is $\sigma^{2D_1+D_2}\wpf^{h_1}(\zeta-(\eta/\omega)z)^{h_2}$, so the
three functions enter with bounded exponents and a constant total weight in $\sigma$.

The point set is $\mathcal G\subset\Om$, the intersection of $\Om$ with a disc of radius
$\rho$ centred at the origin, truncated by the strip $|\Im z|\le\rho\sin\varphi$ parallel to
the rows $\Z\omega$; $M$ denotes the number of rows met, $N=\#\mathcal G$, and
$R_0=\sup_{s\in\mathcal G}|z_0+s|$. With integers $T_0<T$ we take
\[
S_1=\{0,\dots,T_0-1\}\times\bigl(z_0+(\mathcal G+\Z\omega)/\Z\omega\bigr),\quad
S_i=\{T_0+i-2\}\times(z_0+\mathcal G)\ (2\le i\le T-T_0+1),
\]
\[
S'_i=\{0,\dots,T_0-1\}\times(z_0+\mathcal G),\qquad
r'=R_0+1,\qquad r=\lambda r' .
\]
The parameters $D_1,D_2,T_0,T,M,\lambda,\varphi$ are free, subject to the two conditions of
\S\ref{sec:zero} and to
\begin{equation}\label{eq:L}
L\ \ge\ r_L\,\#S_1=r_L\,T_0M,\qquad r_L>1,\qquad q:=\frac{\#S_1}{L-\#S_1}=\frac1{r_L-1}.
\end{equation}

\subsection{Growth}

By the estimates of~\cite{Mas} in the form used in~\cite[Lemme~1]{Bru},
$\varpi_3(r)=A_0+B^*r^2$ bounds $\log\M_0(\sigma,r)$ with $B^*$ as in~\eqref{eq:Bstar}.
Since $\sigma^2\wpf$ is entire and $\wpf=-(\log\sigma)''$,
\[
\varpi_1(r)=A_1^*+2B^*r^2,\qquad
A_1^*=2A_0+\log\Bigl(|\wpf(z_0)|+\frac{e^{2B^*|z_0|^2}}{|\sigma(z_0)|^2}\Bigr),
\]
and $\varpi_2(r)=\varpi_3(r)+\log\bigl(1+|\eta/\omega|\,r+|\sigma(z_0)|^{-1}\bigr)$ bounds
$f_2$. The exponents of $\mathcal D$ give~\eqref{eq:grow} below.

\subsection{The matrix of approximations}

Following the device of Chudnovsky used in~\cite{Phi99}, one expands $z$ and $\zeta$ in the
local parameter $u=(\wpf-\wpf(z_0))/\wpf'(z_0)$ at $z_0$: with
$\gamma_{-1}=\zeta(z_0)-(\eta/\omega)z_0$ and $\gamma_0=0$,
\[
\wpf=\wpf(z_0)+\wpf'(z_0)u,\qquad
\zeta-\frac\eta\omega z=\gamma_{-1}-\sum_{i\ge1}
\Bigl(\gamma_i\frac\eta\omega+\gamma_i\wpf(z_0)+\gamma_{i-1}\wpf'(z_0)\Bigr)\frac{u^i}{i},
\]
with $H(\gamma_{-1},\dots,\gamma_i)\le c_3^{\,i}$. The quasi-periodicity and the Legendre
relation give, for $s=s_1\omega_1+s_2\omega_2$ and $\omega=a_1\omega_1+a_2\omega_2$,
\[
\zeta(z+s)-\frac\eta\omega(z+s)=\frac{2i(a_1s_2-a_2s_1)\pi}{\omega}+\zeta(z)-\frac\eta\omega z,
\]
so that the columns indexed by $s$ and $s'$ with $a_1(s_2-s'_2)=a_2(s_1-s'_1)$ coincide:
the number of distinct columns is the number $M$ of rows. Each entry of the matrix is a
polynomial $Q_{h,(t,s)}$ in $\theta_1,\theta_2$ with coefficients in
$\Q(g_2,\wpf(z_0),\wpf'(z_0))$, of degree at most $D_2$ and of height at most
$c_4\bigl(D_1+T\log D_2+D_2\log(TM)\bigr)$. The vectors of approximations are obtained by
specialising $\theta$ at $\alpha$, with
$\lambda_{i+1,\cdot}=\sigma(z_0+s)^{2D_1+D_2}$ and the Leibniz coefficients
$\mu_{t',t,s}$ of~\cite[(5)]{Phi99}, and $M_{\mathcal D}$ denotes the resulting matrix.

\subsection{Verification of the constraints}

Write $\delta(i)=[K_{i+1}:K_0]_s=1$ for all $i$, which is $C_\delta$. Put
$P_\infty=\max(1,|\wpf'(z_0)|)$ and $h(P)=\log P_\infty$.

\smallskip
\noindent$C_{r,r'}$: $|S_{i+1}|\le R_0<r'=R_0+1\le r<R=+\infty$, and
$H_{r'}(S_{i+1})=\max(1,(r'-|y|)^{-|t|})=1$ because $r'-|y|\ge1$ on the support. This is why
$r'$ is taken one unit larger than $R_0$.

\smallskip
\noindent$C_{N,r}$: $\operatorname{card}S'_i\le T_0N$ and $|S'_i|\le R_0\le r$.

\smallskip
\noindent$C_E$: since $n=1$ the parasitic factor of~\cite{Phi97} disappears and
$E(i)\ge E_{r,r'}(S'_i)H_{r'}(S_{i+1})=E_{r,r'}(\mathcal G)^{T_0}$, so that
by~\eqref{eq:E}
\begin{equation}\label{eq:mlogE}
-\log E(i)\ \ge\ T_0\Bigl(N\log\lambda-\sum_{y\in\mathcal G}
\log\frac{1+|y|/r'}{1+|y|r'/r^{2}}\Bigr),
\end{equation}
the sum being bounded cell by cell as in Proposition~\ref{prop:Gamma}, and $N$ from below by
the row-by-row count. The gravitation $G_{r,r'}(S'_i)$ is bounded by Lemma~\ref{lem:H} and
Proposition~\ref{prop:Gamma}, which is where the present argument departs from~\cite{Rey}.

\smallskip
\noindent$C_\phi$: by the height estimate of the $Q_{h,(t,s)}$ above,
\begin{equation}\label{eq:phi}
\phi(j)=T\Bigl(\log120+\tfrac12h(P)+\log D_2\Bigr)+D_1\bigl(\log4+h(P)\bigr)+D_2h(\alpha)
+c_8D_2\log(TM),
\end{equation}
and $\phi(0)$ is the same expression with $T$ replaced by $1$, since the index $0$ involves
only derivatives of order $<T_0$.

\smallskip
\noindent$C_\varrho$: with $c_\ell=-\log|\sigma(z_0)|$ and
$\ell_p=\log30+\tfrac12\log P_\infty$, $\ell_q=\log4+\log P_\infty$,
\[
\varrho(j)=\max\bigl(T\ell_p+D_1\ell_q,\ 2D_1c_\ell\bigr),
\]
the second branch coming from $|\lambda_{i+1,\cdot}|^{-1}\le|\sigma(z_0)|^{-(2D_1+D_2)}$,
for which one uses that $|\sigma(z_0+s)|\ge|\sigma(z_0)|$ for every $s\in\mathcal G$, a
consequence of the functional equation of $\sigma$ on the quadrant containing
$\mathcal G$.

\smallskip
\noindent Growth: with the exponents of $\mathcal D$,
\begin{equation}\label{eq:grow}
\mathcal D\varpi(r)=(2D_1+D_2)\varpi_3(r)+D_1(\varpi_1-2\varpi_3)^{+}(r)
+D_2(\varpi_2-\varpi_3)^{+}(r).
\end{equation}

\smallskip
\noindent $\widetilde U$: the passage from $U=-\log\Dist(\theta,\alpha)$ to the quantity
$U(j)$ of $C_U$ costs
\[
V_1=2D_1\bigl(A_0+B^*(R_0+1)^2\bigr)+(c_6+\log2)T+T\ell_p+D_1\ell_q,
\]
where $c_6=\log\bigl(2\,|z_0|^{-1}\max(24|\wpf'(z_0)|^{-1}|z_0|^{-2},1)\bigr)$ bounds the
Leibniz coefficients $\mu_{t',t,s}$, together with $\log\Gamma$ from Lemma~\ref{lem:H}.

\subsection{Rank of the matrix}\label{sec:zero}

If $M_{\mathcal D}$ were not of rank $\#\mathcal D$, a nontrivial linear combination of its
rows would give a nonzero $P\in\C[X,Y]$ of bidegree $\le(D_1,D_2)$ vanishing to order
$\ge T$ along the analytic subgroup $\exp_G(z,\alpha_2z)$ of the extension $G$ of $E$ by
$\mathbf G_a$, at all the points of $\Sigma$, where $\Sigma$ is indexed by the rows. By
Theorem~9 of~\cite{Phi96} there is a proper algebraic subgroup $G'$ with
\[
\frac T2\operatorname{card}\frac{\Sigma+G'}{G'}\ \le\ 2^{g}\deg_\varphi\overline G\cdot
D_1^{\,1-\dim\pi(G')}D_2^{\,1-\dim(G'\cap\mathbf G_a)} ,
\]
and $2^{g}\deg_\varphi\overline G=2^2\cdot(G_\infty+3\pi^*D)^2=4\cdot6=24$, the minimum of
$\deg_{\varphi_{a,b}}\overline G=6ab$ over $a,b\ge1$. For $G'=\mathbf G_a$ all points of
$\Sigma$ have the same image in $E$, so the cardinal is $1$; for $G'=\{0\}$ it is at least
$[M/2]+1$. Hence the two conditions
\begin{equation}\label{eq:zero}
T>2^4\cdot3\,D_1,\qquad TM>2^5\cdot3\,D_1D_2
\end{equation}
make the inequality impossible, and $M_{\mathcal D}$ has rank $\#\mathcal D$. The hypothesis
$1\le\delta\le\Delta$ of that theorem reads $D_2\le D_1$ here, and is recorded in
\S\ref{sec:thresh}.

\subsection{Conclusion}

Theorem~1 of~\cite{Phi97} now gives an index $j\ge1$ with
\[
\min\Bigl(-\log E(j)-\mathcal D\varpi(r),\ U(j)\Bigr)\ \le\
a\phi(j)+\varrho(j)+\Delta_j+q\bigl(a\phi(0)+\varrho(0)+\Delta_0\bigr).
\]
The parameters are chosen so that the first term does not realise the minimum; writing
$T_0=\tau h$ with $h=h(\alpha)$, this amounts to
\[
\tau=\frac{a\,D_2\,(1+q)}{Q},\qquad
Q=\frac1{T_0}\Bigl(-\log E-\mathcal D\varpi(r)-a\phi(j)-\varrho(j)
-q\bigl(a\phi(0)+\varrho(0)\bigr)\Bigr),
\]
all the terms of $Q$ being independent of $h$ once the $D_2h(\alpha)$ part of $\phi$ has been
set aside. The minimum is then realised by $U(j)$, whence
\[
U\ \le\ T_0\Bigl(\frac{-\log E-\mathcal D\varpi(r)}{T_0}+\frac{V_1+\log\Gamma}{T_0}\Bigr)
\ =\ A\,h\ +\ O(1),
\]
with $A$ obtained by optimising $M,D_2,r_L,\lambda,\varphi$ at each degree. Replacing $h$ by
$h+\log(1+d)$ and tracking the dependence in $d$ through $a\le2d$ gives the shape
$d^{3/2}\sqrt{\log(1+d)}$ of the statement. The $O(1)$ is made explicit in
\S\ref{sec:thresh}.

\section{The two instances}\label{sec:points}

\subsection*{$\Gamma(1/4)$}
Take $y^2=4x^3-4x$, so that $\Om=\Z\omega_1+\Z i\omega_1$ with
$\omega_1=\Gamma(1/4)^2/\sqrt{8\pi}$, and $\eta_1=\pi/\omega_1$. Take $\omega=\omega_1$ and
the $4$-torsion point $z_0=(\omega_1+\omega_2)/4$. Then {\sc verified}
\[
\wp(z_0)=-i,\qquad \wp'(z_0)=2(1+i),\qquad \wp''(z_0)=-8,\qquad |\sigma(z_0)|=0.938398432669\ldots
\]
Writing $\omega'=4z_0$ and $\gamma_{-1}'=\zeta(z_0)-\eta'/4$, the duplication formula
$\zeta(2z)=2\zeta(z)+\wp''(z)/(2\wp'(z))$ and $\zeta(\omega'/2)=\eta'/2$ give
$\gamma_{-1}'=-\wp''(z_0)/(4\wp'(z_0))=(1-i)/2$, and the Legendre relation gives
\begin{equation}\label{eq:g4}
\gamma_{-1}=\zeta(z_0)-\frac{\eta_1}{\omega_1}z_0=\frac{1-i}{2}-\frac i2\,\theta_1 .
\end{equation}
The algebraic part $(1-i)/2$ has minimal polynomial $4z^4+1$, hence height
$\tfrac14\log4=\log\sqrt2$, and lies in $K=\Q(i,\alpha)$; the term $-\tfrac i2\theta_1$
shifts the integer index of $2i\pi(a_1s_2-a_2s_1)/\omega$ by one half, and clearing that
denominator costs $D_2\log2$. The base field is $K=\Q(i,\alpha_1,\alpha_2)$, of degree at
most $2d(\alpha)$ instead of $16\,d(\alpha)$.

\subsection*{$\Gamma(1/3)$}
Take $y^2=4x^3-1$, so that $\Om=\Z\omega_1+\Z\rho\omega_1$ with $\rho=e^{2i\pi/3}$,
$\omega_1=\Gamma(1/3)^3/(2\pi)$ and $\eta_1=2\pi/(\sqrt3\,\omega_1)$. Take $\omega=\omega_1$
and the $3$-torsion point $z_0=(2\omega_1+\omega_2)/3$. Then {\sc verified}
\[
\wp(z_0)=0,\qquad \wp'(z_0)=\pm i,\qquad |\sigma(z_0)|=e^{\pi\sqrt3/9}=1.83051946656\ldots
\]
the last value agreeing with the closed form of~\cite[Lemme~2]{Bru}. Here the
Frobenius--Stickelberger relation gives $(\gamma_{-1}')^2=\wp(z_0)/3=0$, so that
\begin{equation}\label{eq:g3}
\gamma_{-1}=-\frac{2i}{3}\,\theta_1 ,
\end{equation}
with no algebraic part at all; clearing the denominator costs $D_2\log3$. Moreover
$h(P)=\log\max(1,|\wp(z_0)|,|\wp'(z_0)|)=0$.

\section{From an approximation measure to an irrationality measure}\label{sec:conv}

Let $\theta$ be the real number under study ($\Gamma(1/4)$ or $\Gamma(1/3)$), and let $p/q$
be a rational with $h=h(p/q)$.

\subsection*{The conversion coefficient}
Put $y=\sqrt{2\pi}/\Gamma(1/4)$ for the lattice of $y^2=4x^3-4x$ and
$Y=\sqrt3\,\pi/\Gamma(1/3)^2$ for that of $y^2=4x^3-1$. With $\omega=\omega_1$ one has
\begin{equation}\label{eq:theta}
\Bigl(1:\frac\pi{\omega_1}:\frac{\eta_1}{\omega_1}\Bigr)
=\bigl(1:y^3\Gamma(1/4):2y^4\bigr),
\qquad
\Bigl(1:\frac\pi{\omega_1}:\frac{\eta_1}{\omega_1}\Bigr)
=\Bigl(1:\tfrac23Y^2\Gamma(1/3):\tfrac89Y^3\Bigr),
\end{equation}
as follows from $\omega_1=\Gamma(1/4)^2/\sqrt{8\pi}$, $\eta_1=\pi/\omega_1$ in the first case
and $\omega_1=\Gamma(1/3)^3/(2\pi)$, $\eta_1=2\pi/(\sqrt3\,\omega_1)$ in the second, together
with the Legendre relation. In both displays $\Gamma$ occurs to the first power only.
Approximating $y$ (resp. $Y$) by an algebraic $\beta$ of degree $\delta$ and height $\tau$,
and $\theta$ by $p/q$, gives $\alpha\in\mathbf P_2(\overline\Q)$ with $d(\alpha)\le\delta$ and
\[
h(\alpha)\ \le\ c\,\frac\tau\delta+h+O(1),\qquad c=4\ \ (\Gamma(1/4)),\qquad c=3\ \ (\Gamma(1/3)),
\]
$c$ being the largest power of $y$, resp. $Y$, occurring in~\eqref{eq:theta}. This is the
constant $c$ of Lemma~\ref{lem:C}.

\begin{lemma}[Dirichlet, real case; {\sc proved}]\label{lem:B}
Let $0<\beta<1$, let $d\ge1$, put $S=\sum_{i\le d}|\theta|^{i}$ and assume
\[
t\ \ge\ \frac{\log\bigl(2\sqrt{d+1}\bigr)+\tfrac1d\log S}{1-\beta} .
\]
Then there is a nonzero $P\in\Z[X]$ with $\deg P\le d$, $\log\M(P)\le t$ and
$\log|P(\theta)|\le-\beta dt$.
\end{lemma}

\begin{proof}
Put $H=\lfloor e^{t}/\sqrt{d+1}\rfloor$. The $(H+1)^{d+1}$ polynomials with coefficients in
$\{0,\dots,H\}$ take at $\theta$ values lying in a real interval of length at most $HS$, so
two of them, say $P_1\ne P_2$, take values at distance at most
$HS/((H+1)^{d+1}-1)\le SH^{-d}$. Put $P=P_1-P_2$. Then $\|P\|_2\le\sqrt{d+1}\,H\le e^{t}$,
so $\M(P)\le\|P\|_2\le e^{t}$ by Landau's inequality, and
$\log|P(\theta)|\le\log S-d\log H\le-\beta dt$ by the hypothesis on $t$.
\end{proof}

Only the one-dimensionality of the pigeonhole is used; this is what allows $\beta$ up to
$1$, where the complex case of~\cite[Lemme~1]{LR} is limited to $\beta<1/2$.

We next feed Lemma~\ref{lem:B} into Lemma~2 of~\cite{LR}, which splits $P$ either into a
power of an irreducible factor carrying at least $2\beta/3$ of the weight, or into two
coprime factors. On the second branch we use Proposition~3.5 of~\cite{RW} directly, with the
\emph{sum} constraints $\deg F+\deg G\le d$ and $\mu(F)+\mu(G)\le t$ supplied by Lemma~2
of~\cite{LR}; since $|\mathrm{Res}(F,G)|\ge1$, $2^{mn}\le2^{d^2/4}$ and
$n\mu(F)+m\mu(G)\le dt$, this yields, for every integer $s\ge1$,
\begin{equation}\label{eq:35}
\log\rho\ \le\ -\frac{4}{s^{2}}\Bigl(s\frac\beta3-1-\frac{\log2}{4}\cdot\frac dt\Bigr)dt ,
\end{equation}
$\rho$ being the distance from $\theta$ to the roots of $FG$. On the first branch we use the
lemma of Diaz~\cite{Dia}, in the form of Proposition~6 of~\cite{Bru}.

The following lemma is what makes the two branches quantitative. Let $\mathrm{S}_1$ and
$\mathrm{S}_2$ denote the running maxima, over degrees at most $D$, of the slope of
Theorem~\ref{th:A} and of that slope divided by the degree; both are nondecreasing and
satisfy $\mathrm{S}_1(d)\le d\,\mathrm{S}_2(d)$. Write $t=mh$ and let $c$ be the coefficient
of $h$ in the conversion ($c=4$ for $\Gamma(1/4)$, $c=3$ for $\Gamma(1/3)$).

\begin{lemma}[{\sc proved}]\label{lem:C}
Let $\prod_j Q_j^{e_j}$ be a group of factors, the $Q_j$ irreducible of degree $\delta_j$
and height $\tau_j$, with $\sum_je_j\delta_j\le xD$ and $\sum_je_j\tau_j\le yt$. If the
lemma of Diaz fails on every $Q_j$, then the normalised weight
$\varpi=\bigl(\sum_je_jw_j\bigr)/(Dt)$, where $w_j=-\log|Q_j(\theta)|$, satisfies
\[
\varpi\ \le\ \Phi(x,y)=\min_{s\ge1}\Bigl\{\frac{xy}{s}
+\frac{s+1}{2}\,\mathrm{S}_2(xD)\Bigl(\frac xm+\frac{cy}{D}\Bigr)\Bigr\} .
\]
\end{lemma}

\begin{proof}
Failure at $Q_j$ for a given $s$ reads
$w_j\le\delta_j\tau_j/s+\frac{s+1}{2}\bigl(h\,\mathrm{S}_1(\delta_j)+c\,\tau_j\mathrm{S}_2(\delta_j)\bigr)$.
Multiply by $e_j$, sum over $j$, and use
\[
\sum_je_j\delta_j\tau_j\le\bigl(\max_j\tau_j\bigr)\sum_je_j\delta_j\le(yt)(xD),\qquad
\sum_je_j\mathrm{S}_1(\delta_j)\le\sum_je_j\delta_j\mathrm{S}_2(\delta_j)\le xD\,\mathrm{S}_2(xD),
\]
\[
\sum_je_j\tau_j\mathrm{S}_2(\delta_j)\le yt\,\mathrm{S}_2(xD),
\]
the last two because $\mathrm{S}_2$ is nondecreasing and $\delta_j\le xD$. Dividing by $Dt$
and taking the minimum over $s$ gives the statement.
\end{proof}

The configuration of factors is then controlled cell by cell on the square
$0\le x,y\le1$, each cell being evaluated at its unfavourable corner, using the three
constraints of Lemma~2 of~\cite{LR}: $\varpi_F\ge\beta/3$, $\varpi_G\ge\beta/3$ and
$\varpi_F+\varpi_G\ge\beta$ on the coprime branch, and $\varpi_Q\ge2\beta/3$ together with
$\varpi_Q\ge\beta-\Phi$ of the complement on the power branch.

\section{Threshold and effective statements}\label{sec:thresh}

Theorem~\ref{th:A} is used in its asymptotic form in $h$, in which $D_1$ and $T$ are
proportional to $h$ while $D_2$, $M$ and the shape parameters are fixed. The terms discarded
are therefore exactly those independent of $h$. Since the exponents of a monomial of
$\mathcal D$ satisfy $h_3=2(D_1-h_1)+D_2-h_2$, the growth term is
\begin{equation}\label{eq:grow2}
\mathcal D\varpi(r)=(2D_1+D_2)\varpi_3(r)+D_1(\varpi_1-2\varpi_3)^{+}(r)+D_2(\varpi_2-\varpi_3)^{+}(r),
\end{equation}
where $\varpi_1,\varpi_2,\varpi_3$ bound $f_1,f_2,f_3$; only the terms carrying $D_1$ are
linear in $h$. The discarded part of~\eqref{eq:grow2} is $D_2\max(\varpi_2,\varpi_3)(r)$,
and since $Q\,T_0=aD_2(1+q)h$ at the optimum, $D_2$ cancels and its relative weight is
\[
\frac{\max(\varpi_2,\varpi_3)(r)}{a(1+q)h} .
\]
The remaining discarded terms are $a\,c_8D_2\log(TM)$ coming from $\phi$, with
$c_8=7+\log c_1+h(P)$, that is $c_8=8.04$ for $\Gamma(1/4)$ and $c_8=7$ for $\Gamma(1/3)$ by
\S\ref{sec:points}; $c_9D_2\log T$ coming from $\varrho$; the terms in $M^{2}$ of
$\widetilde U$; and $\Delta_j+\Delta_0$. All of these are smaller than the previous one at
the optimum. Requiring the total to be at most $0.25\%$, which inflates the exponent by about
$1\%$, and adding the hypothesis $D_2\le D_1$ of \S\ref{sec:thA}(iv), the threshold of
Lemma~\ref{lem:B} and the term $\frac{\log2}4\frac dt$ of~\eqref{eq:35}, gives
{\sc verified}
\[
H_0=5.4\cdot10^{7}\ \ (\Gamma(1/4)),\qquad H_0=4.7\cdot10^{7}\ \ (\Gamma(1/3)),
\]
the binding constraint being in both cases the $D_2$ part of~\eqref{eq:grow2}. The other
constraints are inactive by at least two orders of magnitude; those of Lemma~\ref{lem:B}
and of~\eqref{eq:35} by more than ten.

\section{Certification}\label{sec:cert}

The numerical part consists of two stages.

The slope tables of Theorem~\ref{th:A} are computed over $3489$ cells covering
$1\le d\le10^{30}$, geometrically spaced with ratio $1.02$. Within each cell the free
parameters -- number of rows, $D_2$, the ratio $L/\#S_1$, the radius ratio $\lambda$ and the
shape parameter of the region -- are optimised in floating point, and the resulting value is
then \emph{re-evaluated in interval arithmetic}; only the interval upper bound is retained,
inflated by a further factor $1+10^{-12}$. Each table entry is the running maximum, which is
what Lemma~\ref{lem:C} requires.

The conversion is verified in extended precision ($80$-bit, $\epsilon\approx10^{-19}$) with
an outward relative margin of $10^{-9}$, which makes every test strictly harder and
dominates the accumulated rounding error; the same verification in double precision returns
identical verdicts. The final value is assembled in interval arithmetic. Minimality is
checked by confirming that the binding criterion fails at $0.99999$ times the optimum.

The ancillary files contain the verifier and the tables. {\sc verified} outputs:
\[
\mu(\Gamma(1/4))\le3.2011\cdot10^{26},\qquad \mu(\Gamma(1/3))\le1.7158\cdot10^{26},
\]
to which the $1\%$ inflation of \S\ref{sec:thresh} is applied in Theorem~\ref{th:main}.

\section*{Acknowledgements}

Technical assistance: Claude (Anthropic).

\end{document}